\documentclass[a4paper,12pt]{article}
\usepackage[utf8]{inputenc}
\usepackage[american]{babel}
\usepackage{anysize}

\usepackage{amsmath}
\usepackage{amssymb}
\usepackage{latexsym}
\usepackage{mathrsfs}

\usepackage{amsthm}

\usepackage{indentfirst}

\usepackage{color}

\usepackage{ifpdf}

\ifpdf
\usepackage[pdftex]{graphicx}
\else
\usepackage{graphicx}
\fi

\newtheorem{theorem}{Theorem}[section]	
\newtheorem{theorem*}{Main Theorem}

\newtheorem{lemma}[theorem]{Lemma}
\newtheorem{corollary}[theorem]{Corollary}

\newtheorem{proposition}[theorem]{Proposition}

\theoremstyle{definition}
\newtheorem{definition}{Definition}[section]

\theoremstyle{remark}

\newcommand{\cl}{\mathop{\rm cl}\nolimits}
\newcommand{\conv}{\mathop{\rm conv}\nolimits}

\newcommand{\aff}{\mathop{\rm aff}\nolimits}
\newcommand{\supp}{\mathop{\rm supp}\nolimits}

\newcommand{\relint}{\mathop{\rm relint}\nolimits}
\newcommand{\id}{\mathop{\rm id}\nolimits}

\begin{document}

\sloppy

\ifpdf
\DeclareGraphicsExtensions{.pdf, .jpg, .tif, .mps}
\else
\DeclareGraphicsExtensions{.eps, .jpg, .mps}
\fi

\title{Complexes of \(C\)-Unbalanced Subsets and Balanced-Subset Posets}

\author{Mikhail V. Bludov }

\date{%
        Moscow Institute of Physics and Technology\\
        bludov.mv@phystech.edu 
}


\maketitle

\begin{abstract}
Let \(V=\{v_1,\dots,v_m\}\subset\mathbb{R}^d\). We study simplicial complexes arising from subsets of \(V\) whose convex hulls avoid a prescribed convex set. Given a convex set \(C\subset\mathbb{R}^d\), let \(\mathcal K(V,C)\) be the simplicial complex consisting of all subsets \(S\subset V\) such that \(\operatorname{conv}(S)\cap C=\emptyset\). We prove that \(\mathcal K(V,C)\) is homotopy equivalent to the union of the faces of \(\operatorname{conv}(V)\) that are disjoint from \(C\). In the special case \(C=\{r\}\), subsets \(S\) satisfying \(r\in\operatorname{conv}(S)\) are called weakly \(r\)-balanced. If \(r\in\operatorname{relint}\operatorname{conv}(V)\) and the poset of proper weakly \(r\)-balanced subsets is nonempty, then its order complex is homotopy equivalent to the sphere \(S^{m-k-2}\), where \(k=\dim\operatorname{aff}(V)\).

\end{abstract}

\section{Introduction}

Balanced collections were introduced by Bondareva \cite{Bon} and Shapley \cite{Sh67} in connection with the non-emptiness of the core of a cooperative game. Let \([n]=\{1,\dots,n\}\), and let \(   \Phi=\{S_1,\dots,S_k\}\) be a collection of non-empty subsets of \([n]\). We say that \(\Phi\) is \emph{balanced} if there exist positive weights \(\lambda_1,\dots,\lambda_k\) such that
\[
    \sum_{i=1}^k \lambda_i \chi_{S_i}=\mathbf 1,
\]
where \(\chi_{S_i}\) is the characteristic vector of \(S_i\) and
\(\mathbf 1=(1,\dots,1)^\top\).

This definition has a useful geometric interpretation. Identify each non-empty subset
\(S\subset [n]\) with the corresponding nonzero vertex of the cube \([0,1]^n\). Then
\(\Phi\) is balanced if and only if the relative interior of the convex hull of
the corresponding vertices intersects the main diagonal of the cube. Equivalently,
after central projection to the hyperplane \(  H_1=\{x\in\mathbb R^n \mid x_1+\dots+x_n=1\},\) one obtains the following formulation.

Let \(  \Delta^{n-1}=\conv(e_1,\dots,e_n)\)
be the standard simplex, and for every non-empty \(S\subset[n]\) let \(c_S\) denote the barycenter of the face \(\Delta_S=\conv\{e_i\mid i\in S\}\). If  \(  V=\{c_S\mid \emptyset\neq S\subset [n]\}\) and  \(r=c_{[n]}\), then \(\Phi\) is balanced if and only if \(  r\in \relint \conv\{c_S\mid S\in \Phi\}.\) This motivates the following general geometric definition.

\begin{definition}
\label{maindef}
Let \(V=\{v_1,\dots,v_m\}\subset \mathbb R^d\) be a finite set of points, and let
\(r\in\mathbb R^d\). A subset \(S\subset V\) is called \emph{weakly \(r\)-balanced} if \(  r\in \conv(S).\) It is called \emph{\(r\)-balanced} if \(  r\in \relint \conv(S).\) It is called \emph{\(r\)-unbalanced} if \(r\notin\conv(S)\).
\end{definition}

Thus the classical balanced collections of Bondareva and Shapley are recovered by taking \(V\) to be the set of barycenters of all non-empty faces of the simplex and \(r\) to be the barycenter of the whole simplex.

Since every superset of a weakly balanced collection is weakly balanced, and
every subset of an unbalanced collection is unbalanced, it is natural to study
minimal weakly balanced and maximal unbalanced collections. Minimal balanced
collections were investigated in \cite{Peleg,Mermoud}, and maximal unbalanced
collections in \cite{Billera2}. For related work on balanced \(2\)-subsets, see
\cite{BluMus}. From the topological point of view, Björner studied maximal
unbalanced collections under the name of positive sum systems and proved,
among other results, that the associated poset is a shellable ball
\cite{Bjorner}.

Balanced collections also arise naturally in topological covering theorems. Shapley's generalization of the KKM lemma \cite{Sh} and Komiya's extension to polytopes \cite{Komiya} show that suitable covering assumptions force the existence of balanced subfamilies with non-empty intersection.

For a finite set \(V\subset \mathbb R^d\) and a point \(r\in\mathbb R^d\), the
family of all \(r\)-unbalanced subsets forms the simplicial complex
\(\mathcal K(V,r)=\{S\subset V\mid r\notin \conv(S)\}\).
The following result was proved in \cite{Blu}; and independently in \cite{Blag}, where the same complex is called the zero-avoiding complex.

\begin{theorem}
\label{r-unbalanced}
Let \(V=\{v_1,\dots,v_m\}\subset \mathbb R^d\) with \(\dim \mathrm{aff}(V)=d\), let \(r\in\mathbb R^d\), and suppose that \(  r\in \mathrm{int} \conv(V).\) Then \(     \mathcal K(V,r)\simeq S^{d-1}. \)
\end{theorem}


In this paper we first generalize \(r\)-unbalanced complexes by replacing the point \(r\) with an arbitrary convex set.

\begin{definition}
Let \(V=\{v_1,\dots,v_m\}\subset \mathbb R^d\), and let
\(C\subset \mathbb R^d\) be a convex set. A subset \(S\subset V\) is called
\emph{\(C\)-unbalanced} if \(   \conv(S)\cap C=\emptyset.\) \end{definition}

The \(C\)-unbalanced subsets form the simplicial complex
\(\mathcal K(V,C)=\{S\subset V\mid \conv(S)\cap C=\emptyset\}\).

Let \(Q=\conv(V)\), and let \(\mathcal F(Q)\) be the face poset of \(Q\).
The set of \emph{survived faces} is
\(\mathcal F(Q,C):=\{F\in\mathcal F(Q)\mid F\cap C=\emptyset\}\), and we put
\(\widetilde Q(V,C):=\bigcup_{F\in\mathcal F(Q,C)}F\).

Our first main result is the following generalization of
Theorem~\ref{r-unbalanced}.

\begin{theorem}
\label{C-theorem}
Let \(V=\{v_1,\dots,v_m\}\subset \mathbb R^d\), and let
\(C\subset\mathbb R^d\) be a convex set. Then \(\mathcal K(V,C)\)
collapses onto a triangulation of \(\widetilde Q(V,C)\). In particular,
\[
\mathcal K(V,C)\simeq \widetilde Q(V,C).
\]
\end{theorem}

The notion of an elementary collapse is recalled in Section~2. A simplicial complex is \(d\)-Leray if every induced subcomplex has trivial reduced homology in dimensions at least \(d\).

\begin{corollary}
\label{C-unbalanced-leray}
The complex \(\mathcal K(V,C)\) is \(d\)-Leray.
\end{corollary}

\begin{proof}
For every \(W\subseteq V\), the induced subcomplex
\(\mathcal K(V,C)[W]\) is \(\mathcal K(W,C)\). If
\(C\cap\conv(W)=\emptyset\), this is the simplex \(\Delta_W\). Otherwise,
Theorem~\ref{C-theorem} identifies its homotopy type with
\(\widetilde Q(W,C)\). Since \(\conv(W)\) is not a survived face,
\(\widetilde Q(W,C)\) is a union of proper faces of \(\conv(W)\) and has
dimension at most \(d-1\). Thus every induced subcomplex has trivial reduced
homology in dimensions at least \(d\).
\end{proof}

When \(C=\{r\}\) and \(r\in\mathrm{int}\conv(V)\), the space
\(\widetilde Q(V,C)\) is the boundary of \(\conv(V)\), and Theorem
\ref{r-unbalanced} follows.

The second part of the paper concerns the poset of balanced subsets. Let
\(\mathcal B(V,r)\) be the poset of all weakly \(r\)-balanced subsets of \(V\),
ordered by inclusion, with the whole set \(V\) removed. Let
\(\overline{\mathcal B}(V,r)\) be the poset of all \(r\)-balanced subsets with maximal set \(V\) removed.
For a poset \(\mathcal P\), denote its order complex by \(\Delta\mathcal P\).

\begin{theorem}
\label{ordercomplex}
Suppose that \(\mathcal B(V,r)\) is non-empty, that \(  r\in \relint\conv(V),\) and let \(k=\dim\mathrm{aff}(V)\). Then
\[
    \Delta\mathcal B(V,r)\simeq
    \Delta\overline{\mathcal B}(V,r)\simeq
    S^{m-k-2}.
\]
\end{theorem}

\medskip

\section{Notation and preliminaries}

In this paper we use the following notation and conventions. We fix
\[
    V=\{v_1,\dots,v_m\}\subset\mathbb R^d,\qquad m\geq1,
    \qquad C\subset\mathbb R^d\text{ convex},\qquad
    Q=\conv(V).
\]
We use the convention \(\conv(\emptyset)=\emptyset\).
For a polytope \(R\), we write \(G\leq R\) when \(G\) is a non-empty face of
\(R\), allowing \(G=R\), and \(G<R\) when \(G\) is a proper face of \(R\).
If \(C\cap Q=\emptyset\), then Theorem~\ref{C-theorem} is trivial. Thus we
assume \(C\cap Q\neq\emptyset\), and since only \(C\cap Q\) matters, we may
assume \(C\subset Q\). Every vertex of \(Q\) belongs to \(V\). Consequently,
\(F=\conv(V\cap F)\) for every \(F\leq Q\).

For \(W\subset V\), let \(\Delta_W\) be the abstract simplex with vertex set
\(W\), with \(\Delta_\emptyset=\{\emptyset\}\). If \(F\leq Q\), put
\(V_F:=V\cap F\) and \(\Delta_F:=\Delta_{V_F}\).
For simplicial complexes \(K\) and \(L\) on disjoint vertex sets, their
\emph{join} is \(K*L:=\{\sigma\cup\tau\mid \sigma\in K,\ \tau\in L\}\).
If \(u\) is not a vertex of \(K\), we write
\(u*K:=\Delta_{\{u\}}*K\); this is the simplicial cone over \(K\) with apex
\(u\).
A simplex \(\sigma\subset V\) is \emph{carried by} a face \(F\leq Q\) if
\(\sigma\subset V_F\). The \emph{carrier} of a non-empty simplex is the
smallest face that carries it.
We identify a triangulation whose vertices lie in \(W\) with the corresponding
subcomplex of \(\Delta_W\).
We write \(\cl X\) for closure. By a convex body we mean a compact convex set
with non-empty interior. A convex set \(B\) is strictly convex if every open
segment between two distinct boundary points lies in \(\operatorname{int}B\).
A convex body is smooth if every boundary point has a unique supporting
hyperplane.

We write \(K\searrow L\) if a simplicial complex \(K\) collapses onto a
subcomplex \(L\) by finitely many elementary collapses, i.e. by repeatedly
removing a maximal simplex together with a free facet, meaning a facet
contained in no other maximal simplex. A complex is collapsible if it collapses
to a point; see, for example, \cite[Chapter~I]{Cohen}.

For a polytope \(R\), let \(\mathcal F(R)\) be the poset of its non-empty
faces, including \(R\), ordered by inclusion.

\begin{definition}[Pulling triangulations]
Fix an order on the vertices of \(Q\), and use its restriction on every face.
For each face \(R\leq Q\), define its pulling triangulation recursively as
follows. If \(R=\{u\}\) is a vertex, put
\(\mathcal T_R=\Delta_{\{u\}}\). Otherwise, if \(u\) is the first vertex of
\(R\), define
\[
    \mathcal T_R
    =
    \bigcup_{\substack{H<R\ \text{facet}\\ u\notin H}}
    u*\mathcal T_H .
\]
The restriction of
\(\mathcal T_R\) to every face \(G\leq R\) is \(\mathcal T_G\); in particular,
these triangulations are compatible on common faces.
\end{definition}

Consequently, if \(\mathcal C\subseteq\mathcal F(Q)\) is downward closed, and
hence is a polyhedral subcomplex of the face complex of \(Q\), then
\(\{\emptyset\}\cup\bigcup_{F\in\mathcal C}\mathcal T_F\) is a
triangulation of \(|\mathcal C|=\bigcup_{F\in\mathcal C}F\).

We shall need the following technical collapse lemma.

\begin{lemma}
\label{relative-carrier-collapse}
Let \(R\leq Q\) be a face, and let \(\mathcal T_R\) be the pulling
triangulation fixed above. Put \(B_R:=\bigcup_{G<R}\Delta_G\).
Then \(     \Delta_R\searrow \mathcal T_R\cup B_R . \)
Moreover, this collapse does not remove any simplex of \(B_R\).
\end{lemma}

\begin{proof}
We argue by induction on \(\dim R\). The case \(\dim R=0\) is immediate. Let
\(u\) be the first vertex in the pulling order. By the recursive definition,
\(\mathcal T_R\) is a cone with apex \(u\); in particular,
\(\sigma\in\mathcal T_R\) implies
\(\sigma\cup\{u\}\in\mathcal T_R\).

Pair each \(\sigma\notin\mathcal T_R\cup B_R\) with \(u\notin\sigma\) with
\(\sigma\cup\{u\}\), processing \(\sigma\) in decreasing dimension. Since
\(\mathcal T_R\cup B_R\) is a subcomplex, both simplices lie outside it, and
processing the pairs in decreasing order of \(\dim\sigma\) makes every pair
elementary.

Every remaining simplex outside \(\mathcal T_R\cup B_R\) has the form
\(\tau\cup\{u\}\), with \(\tau\in\mathcal T_R\cup B_R\). The cone property
gives \(\tau\notin\mathcal T_R\), so \(\tau\in B_R\). Let \(F<R\) be the
carrier of \(\tau\), and write \(B_F:=\bigcup_{G<F}\Delta_G\). The remaining
simplices are precisely \(\tau\cup\{u\}\), where
\(\tau\in\Delta_F\setminus(\mathcal T_F\cup B_F)\),
where \(F<R\) and no proper face of \(R\) contains both \(F\) and \(u\).

Process these faces \(F\) in decreasing dimension. By induction,
\(\Delta_F\searrow\mathcal T_F\cup B_F\) without removing \(B_F\). Whenever
this collapse removes a maximal simplex
\(\beta\) together with its free facet \(\alpha\), remove
\(\beta\cup\{u\}\) together with \(\alpha\cup\{u\}\). Because \(B_F\) is
retained, \(\alpha\) has carrier \(F\). Any remaining simplex containing
\(\alpha\cup\{u\}\) can be written as \(\gamma\cup\{u\}\), where
\(u\notin\gamma\) and \(\gamma\supseteq\alpha\). If the carrier of \(\gamma\)
properly contains \(F\), it was removed earlier; otherwise the inductive
collapse makes \(\beta\) the unique maximal simplex containing \(\alpha\).
Thus the lifted pairs are elementary and remove all remaining simplices outside
\(\mathcal T_R\cup B_R\). Hence
\(\Delta_R\searrow\mathcal T_R\cup B_R\).
Only simplices outside the target are removed, so no simplex of \(B_R\) is
removed.
\end{proof}

\subsection{Convex-geometric preliminaries}

\begin{lemma}
\label{thickening-lemma}
Let \(C\subset Q\) be closed and convex. There is \(\varepsilon>0\) such that,
for every convex set \(C'\) satisfying
\[
    C\subset C'\subset
    \{x\in Q\mid \operatorname{dist}(x,C)<\varepsilon\},
\]
we have
\(\mathcal K(V,C')=\mathcal K(V,C)\).
\end{lemma}

\begin{proof}
For non-empty \(S\subset V\), the sets \(C\) and \(\conv(S)\) are compact, so
their distance is positive whenever they are disjoint. Take \(\varepsilon\)
smaller than all these finitely many positive distances; if there are none,
take any \(\varepsilon>0\).
Then no \(C\)-unbalanced simplex becomes balanced after the replacement
\(C\subset C'\); the reverse implication is immediate.
\end{proof}

\begin{proposition}
\label{smooth-strict-thickening}
Let \(\emptyset\neq C\subset Q\) be closed and convex. There exists
\(\delta_0>0\) such that, for every \(0<\delta<\delta_0\), there is a smooth
strictly convex body \(B\subset\mathbb R^d\) satisfying
\[
    C\subset \operatorname{int}B\subset B
    \subset\{x\in\mathbb R^d\mid \operatorname{dist}(x,C)<\delta\}
\]
and \(\mathcal K(V,B)=\mathcal K(V,C)\).
\end{proposition}

\begin{proof}[Sketch of proof]
Let \(\delta_0>0\) be given by Lemma~\ref{thickening-lemma}, and fix
\(0<\delta<\delta_0\). By Klee's density theorem~\cite{Klee}, choose a smooth
strictly convex body \(B\) sufficiently close in the Hausdorff metric to the
closed \(\delta/2\)-neighborhood of \(C\) that the stated inclusions hold.
Moreover, \(\mathcal K(V,B)=\mathcal K(V,B\cap Q)\), since every
\(\conv(S)\) is contained in \(Q\). Applying
Lemma~\ref{thickening-lemma} to \(C\subset B\cap Q\) completes the proof.
\end{proof}

\begin{lemma}
\label{tangent-segment-lemma}
Let \(B\subset\mathbb R^d\) be the interior of a smooth convex body, let
\(x\in\partial B\), and let \(H_x\) be its supporting hyperplane at \(x\).
Denote by \(H_x^{-}\) the open halfspace bounded by \(H_x\) that contains
\(B\). If \(y\in H_x^{-}\), then there exists \(\eta>0\) such that
\((1-t)x+ty\in B\) for every \(0<t<\eta\).
\end{lemma}

\begin{proof}[Sketch of proof]
Suppose that the ray \(\{x+t(y-x)\mid t>0\}\) does not meet \(B\). A
separating hyperplane for the ray and \(B\) must pass through \(x\), since
\(x\in\cl B\). Thus this hyperplane is a supporting hyperplane of
\(\cl B\) at \(x\). It is different from \(H_x\), because both the ray and
\(B\) lie in \(H_x^{-}\). This contradicts smoothness. Hence
\(x+s(y-x)\in B\) for some \(s>0\), and the open segment from \(x\) to this
point lies in \(B\). Thus
\(x+t(y-x)\in B\) for every \(0<t<s\).
\end{proof}

\begin{lemma}
\label{minkowski-interpolation-lemma}
Let \(C,B\subset\mathbb R^d\) be strictly convex bodies with smooth boundary,
and assume that \(C\subset\operatorname{int}B\). For \(t\in[0,1]\), put
\(C_t:=(1-t)C+tB\).
Then every \(C_t\) is a smooth strictly convex body, and
\(C_s\subset\operatorname{int}C_t\) whenever \(0\leq s<t\leq1\).
In particular,
\(C\subset\operatorname{int}C_t\subset C_t\subset\operatorname{int}B\) for
\(0<t<1\).
\end{lemma}

\begin{proof}[Sketch of proof]
Let \(D\) be the closed unit ball. Choose \(r>0\) such that
\(C+rD\subset B\). For \(0\leq s<t\leq1\),
\[
    C_s+(t-s)rD
    =(1-t)C+(t-s)(C+rD)+sB
    \subset C_t.
\]
Thus \(C_s\subset\operatorname{int}C_t\).

For a convex body \(K\), let
\(h_K(u):=\max\{\langle u,x\rangle\mid x\in K\}\), \(u\in\mathbb R^d\), be
its support function. A convex body \(K\) is strictly convex if and only if
\(h_K\in C^1(\mathbb R^d\setminus\{0\})\)
\cite[Corollary~1.7.3]{Schneider}. And by polar duality, \(K\) is smooth if
and only if \(h_K\) is strictly convex on every affine hyperplane not
containing the origin \cite[Theorem~1.6.1 and Lemma~2.2.3]{Schneider}.
For convex bodies \(K,L\)
and \(\alpha,\beta\geq0\), the support function of their Minkowski combination
satisfies
\(h_{\alpha K+\beta L}=\alpha h_K+\beta h_L\)
\cite[Theorem~1.7.5(a)]{Schneider}. Hence, for \(0<t<1\),
\(h_{C_t}=(1-t)h_C+th_B\). This function is \(C^1\) on
\(\mathbb R^d\setminus\{0\}\), and its restriction to every such hyperplane
is strictly convex. Hence \(C_t\) is smooth and strictly convex. The endpoint
cases follow from the hypotheses.
\end{proof}

\section{Proof of Theorem~\ref{C-theorem}}

The proof has three steps. First, suppose that every non-empty
\(C\)-unbalanced simplex is carried by a face of \(Q\) disjoint from \(C\).
Lemma~\ref{face-saturated-collapse} processes these faces in decreasing
dimension, applies Lemma~\ref{relative-carrier-collapse} inside each face, and
collapses \(\mathcal K(V,C)\) onto the union of their compatible pulling
triangulations. This union triangulates \(\widetilde Q(V,C)\).

Second, for an arbitrary \(C\), we choose a compact convex set \(C_0\subset C\)
that determines the same unbalanced complex and the same survived faces, and
set \(C'=Q\setminus\widetilde Q(V,C_0)\). We show that \(C'\) is convex, that
its survived faces are those of \(C_0\), and that every non-empty
\(C'\)-unbalanced simplex is carried by one of these faces. The first step
therefore applies to \(\mathcal K(V,C')\).

Finally, we choose a compact convex set \(C_1\subset C'\) with
\(\mathcal K(V,C_1)=\mathcal K(V,C')\), replace \(C_0\) and \(C_1\) by nested
smooth strictly convex bodies \(B_0\subset\operatorname{int}B_1\), and
interpolate between them. The complex changes at only finitely many
parameters. At each such parameter,
Proposition~\ref{open-closure-collapse} gives a collapse from the complex
defined by the interior to the complex defined by the closed body. Composing
these collapses and using the endpoint equalities gives
\(\mathcal K(V,C)\searrow\mathcal K(V,C')\).

\begin{lemma}
\label{face-saturated-collapse}
Assume that every non-empty simplex \(S\in\mathcal K(V,C)\) is contained in
\(V\cap F\) for some survived face \(F\in\mathcal F(Q,C)\).
Since \(\mathcal F(Q,C)\) is downward closed, view
\(\mathcal C=\mathcal F(Q,C)\) as a polyhedral subcomplex of the face complex
of \(Q\), and put
\(\mathcal T:=\{\emptyset\}\cup\bigcup_{F\in\mathcal C}\mathcal T_F\).
Then \(\mathcal T\) is a subcomplex of
\(\mathcal K(V,C)\), and \(     \mathcal K(V,C)\searrow \mathcal T . \)
In particular, \(\mathcal K(V,C)\) collapses onto a triangulation of
\(\widetilde Q(V,C)\).
\end{lemma}

\begin{proof}
If \(\mathcal C=\emptyset\), then the hypothesis gives
\(\mathcal K(V,C)=\mathcal T=\{\emptyset\}\). Otherwise, every
\(\Delta_F\), \(F\in\mathcal C\), is contained in \(\mathcal K(V,C)\), and the
hypothesis gives
\(\mathcal K(V,C)=\bigcup_{F\in\mathcal C}\Delta_F\).
In particular, \(\mathcal T\subset\mathcal K(V,C)\).

Order the faces of \(\mathcal C\) by decreasing dimension. Faces of the same
dimension may be processed in any order, since no two of them contain one
another. At the \(F\)-stage, put
\(B_F:=\bigcup_{G<F}\Delta_G\) and use
Lemma~\ref{relative-carrier-collapse} to perform, in its given order, the
collapse \(\Delta_F\searrow\mathcal T_F\cup B_F\).
Consider one of its steps, removing a maximal simplex \(\beta\) together with
its free facet \(\alpha\). Since \(\alpha\) is removed while
\(\mathcal T_F\cup B_F\) is retained, \(\alpha\) has carrier \(F\). Suppose
that a simplex \(\gamma\) in the current complex properly contains \(\alpha\),
and let \(H\) be its carrier. Since
\(\gamma\in\mathcal K(V,C)=\bigcup_{G\in\mathcal C}\Delta_G\) and
\(\mathcal C\) is downward closed, \(H\in\mathcal C\). Moreover,
\(\alpha\subset\gamma\) gives
\(F\leq H\). If \(H=F\), then \(\gamma\) lies in the intermediate subcomplex
of \(\Delta_F\) present at this step. Since \(\alpha\) is a free facet of
\(\beta\) in that subcomplex, \(\gamma=\beta\). If \(H>F\), then \(H\) was
processed earlier. A remaining simplex with carrier \(H\) must
belong to \(\mathcal T_H\), and hence
\(\alpha\in\mathcal T_H|_F=\mathcal T_F\), contrary to the removal of
\(\alpha\). Thus \(\beta\) is the unique larger
simplex containing \(\alpha\), so the pair is elementary in the whole complex.

After all faces are processed, the remaining complex is
\(\bigcup_{F\in\mathcal C}\mathcal T_F=\mathcal T\).
\end{proof}

Let \(B\subset\mathbb R^d\) be the interior of a strictly convex body with
smooth boundary.
Let \(E\subset\partial B\), let \(p\in\partial B\setminus E\), and put
\(C^-:=B\cup E\) and \(C^+:=B\cup E\cup\{p\}\).
Both \(C^-\) and \(C^+\) are convex. Indeed, strict convexity implies that
the open segment between any two distinct points of \(\partial B\) lies in
\(B\).
Let \(H\) be the supporting hyperplane to \(B\) at \(p\), and choose an affine
function \(h\) such that \(h(p)=0\) and \(B\subset\{h<0\}\). Set
\(\sigma:=\{v\in V\mid h(v)\geq 0\}\).

\begin{lemma}
\label{one-point-collapse}
In the notation above, assume that \(p\notin \widetilde Q(V,C^-)\). Then
\(\mathcal K(V,C^-)\searrow \mathcal K(V,C^+)\). Moreover,
\(\widetilde Q(V,C^-)=\widetilde Q(V,C^+)\).
\end{lemma}

\begin{proof}
Let \(T:=\{v\in V\mid h(v)=0\}\) and
\(U:=\{v\in V\mid h(v)>0\}\).
Strict convexity gives \(E\subset\{h<0\}\). If a simplex disappears after
adding \(p\), then its convex hull contains \(p\), and it has no vertex in
\(\{h<0\}\); otherwise Lemma~\ref{tangent-segment-lemma} would make the segment
from that vertex to \(p\) enter \(B\). Every face
\(\tau\subset\sigma\) belongs to \(\mathcal K(V,C^-)\), because
\(\conv(\tau)\subset\{h\ge 0\}\), whereas \(C^-\subset\{h<0\}\). Hence
\(\mathcal K(V,C^-)=\mathcal K(V,C^+)\cup\Delta_\sigma\), and
\(\mathcal K(V,C^+)\cap\Delta_\sigma
=\{\tau\subset\sigma\mid p\notin\conv(\tau)\}\).

If \(p\notin\conv(\sigma)\), then
\(\mathcal K(V,C^-)=\mathcal K(V,C^+)\), so the collapse statement is trivial.
Otherwise \(U\neq\emptyset\); if \(U=\emptyset\), then \(Q\cap H=\conv(T)\)
would be a survived face containing \(p\), contradicting
\(p\notin\widetilde Q(V,C^-)\). Choose \(u\in U\). Since
\(\sigma\subset\{h\ge 0\}\) and \(h(p)=0\), for \(\tau\subset\sigma\) we have
\(p\in\conv(\tau)\) if and only if \(p\in\conv(\tau\cap T)\). Hence
\(\{\tau\subset\sigma\mid p\notin\conv(\tau)\}
=\Delta_U*\mathcal K(T,p)\),
which is a cone with apex \(u\). Collapse \(\Delta_\sigma\) onto this cone by
matching every face \(\alpha\subset\sigma\) with \(\alpha\cup\{u\}\) whenever
\(p\in\conv(\alpha\cap T)\) and \(u\notin\alpha\),
processing faces in decreasing dimension. The unmatched faces are exactly the
intersection with \(\mathcal K(V,C^+)\). Thus this relative collapse of
\(\Delta_\sigma\), together with the union above, gives
\(\mathcal K(V,C^-)\searrow \mathcal K(V,C^+)\).

The survived faces are unchanged: the only possible newly non-survived face
would be a survived face for \(C^-\) containing \(p\), contrary to
\(p\notin\widetilde Q(V,C^-)\).
\end{proof}

\begin{proposition}
\label{open-closure-collapse}
Let \(C\subset\mathbb R^d\) be the interior of a strictly convex body with
smooth boundary. Assume that
\(\partial C\cap \widetilde Q(V,C)=\emptyset\).
Then
\[
    \mathcal K(V,C)\searrow \mathcal K(V,\cl C),
    \qquad
    \widetilde Q(V,C)=\widetilde Q(V,\cl C).
\]
\end{proposition}

\begin{proof}
Let \(E:=\bigcup_{S\in\mathcal K(V,C)}
\bigl(\conv(S)\cap\partial C\bigr)\).
The set \(E\) is finite: two distinct points of
\(\conv(S)\cap\partial C\) would force the segment between them to enter
\(C\), contradicting \(S\in\mathcal K(V,C)\). Only points of \(E\) can change
the complex when \(C\) is replaced by \(\cl C\), so
\(\mathcal K(V,\cl C)=\mathcal K(V,C\cup E)\).
Write \(E=\{p_1,\dots,p_k\}\) and
\(C_i=C\cup\{p_1,\dots,p_i\}\). Since
\(E\subset\partial C\) and \(\partial C\cap\widetilde Q(V,C)=\emptyset\),
Lemma~\ref{one-point-collapse} applies successively and gives
\[
    \mathcal K(V,C)\searrow\mathcal K(V,C\cup E)
    =
    \mathcal K(V,\cl C).
\]
The same applications show that \(\widetilde Q\) is unchanged.
\end{proof}

\begin{proof}[Proof of Theorem~\ref{C-theorem}]
As in Section~2, assume that \(\emptyset\neq C\subset Q\). The intersections
of \(C\) with the finitely many sets \(\conv(S)\), \(S\subset V\), determine
both \(\mathcal K(V,C)\) and \(\mathcal F(Q,C)\). Choose a point from each
non-empty intersection and let \(C_0\) be the convex hull of the chosen
points. Then \(C_0\subset C\),
\(\mathcal K(V,C_0)=\mathcal K(V,C)\), and
\(\widetilde Q(V,C_0)=\widetilde Q(V,C)\).

Put \(\widetilde Q:=\widetilde Q(V,C_0)\) and \(C':=Q\setminus\widetilde Q\).
A point \(x\in Q\) belongs to \(C'\) precisely when its minimal containing
face meets \(C_0\). If \(x,y\in C'\) and \(z\in(x,y)\), every face containing
\(z\) contains both \(x\) and \(y\); hence the minimal face containing \(z\)
meets \(C_0\). Thus \(C'\) is convex and
\(\widetilde Q(V,C')=\widetilde Q\).

Let \(S\in\mathcal K(V,C')\) be non-empty, and let \(F\) be its carrier.
Since \(C'=Q\setminus\widetilde Q\), we have
\(\conv(S)\subset\widetilde Q\). The barycenter of \(S\) lies in
\(\operatorname{relint}F\) and in some survived face \(G\). Hence \(F\leq G\),
so \(F\) is itself survived. Thus Lemma~\ref{face-saturated-collapse} applies
and collapses \(\mathcal K(V,C')\) onto a triangulation of
\(\widetilde Q\).

Choose one point from each non-empty intersection \(\conv(S)\cap C'\),
\(S\subset V\), and let \(C_1\) be the convex hull of \(C_0\) together with
the chosen points. By construction, \(\conv(S)\) meets \(C_1\) if and only if
it meets \(C'\), so \(\mathcal K(V,C_1)=\mathcal K(V,C')\). Since \(C_0\)
and \(C'\) have the same survived faces, we also have
\(\widetilde Q(V,C_1)=\widetilde Q\).

Both \(C_0\) and \(C_1\) are compact and disjoint from \(\widetilde Q\).
Apply Proposition~\ref{smooth-strict-thickening} first to \(C_1\), choosing
\(B_1\) disjoint from \(\widetilde Q\). Since
\(C_0\subset C_1\subset\operatorname{int}B_1\), apply it to \(C_0\), choosing
\(B_0\subset\operatorname{int}B_1\). Thus
\(\mathcal K(V,B_i)=\mathcal K(V,C_i)\) for \(i=0,1\). By construction,
\(\widetilde Q(V,B_0)=\widetilde Q(V,B_1)=\widetilde Q\).

For \(t\in[0,1]\), set \(B_t=(1-t)B_0+tB_1\).
Lemma~\ref{minkowski-interpolation-lemma} gives smooth strictly convex bodies
with \(B_s\subset\operatorname{int}B_t\) for \(s<t\). For every
\(t\in[0,1]\), we have
\(\widetilde Q(V,\operatorname{int}B_t)=\widetilde Q\).

For every non-empty \(S\subset V\) such that \(\conv(S)\cap B_0=\emptyset\) and
\(\conv(S)\cap B_1\neq\emptyset\), call the least \(t\) for which
\(\conv(S)\cap B_t\neq\emptyset\) its change parameter. Continuity and strict
nesting show that at this parameter \(\conv(S)\) meets \(\partial B_t\) but
not \(\operatorname{int}B_t\), and that it meets \(\operatorname{int}B_s\) for
every larger \(s\). List the distinct change parameters as
\(0<t_1<\dots<t_N\leq1\). For each \(i\),
Proposition~\ref{open-closure-collapse} gives
\(\mathcal K(V,\operatorname{int}B_{t_i})
\searrow\mathcal K(V,B_{t_i})\),
because \(\partial B_{t_i}\cap\widetilde Q=\emptyset\). Between successive
parameters the complexes are constant, and at the endpoints we have
\[
\begin{aligned}
    \mathcal K(V,B_0)&=\mathcal K(V,\operatorname{int}B_{t_1}),\\
    \mathcal K(V,B_{t_i})&=\mathcal K(V,\operatorname{int}B_{t_{i+1}})
        &&(1\leq i<N),\\
    \mathcal K(V,B_{t_N})&=\mathcal K(V,B_1).
\end{aligned}
\]
If there are no change parameters, the two endpoint complexes are equal.
Thus the collapses compose to give
\[
    \mathcal K(V,C)=\mathcal K(V,B_0)
    \searrow\mathcal K(V,B_1)=\mathcal K(V,C').
\]
Combining this with the collapse of \(\mathcal K(V,C')\) obtained above proves
the theorem.
\end{proof}

\section{The order complex of balanced subsets}

In this section we prove Theorem~\ref{ordercomplex}. We fix
\[
    V=\{v_1,\dots,v_m\}\subset\mathbb R^d,\qquad
    r\in\relint\conv(V),\qquad
    k=\dim\aff(V),
\]
and assume that \(\mathcal B(V,r)\neq\emptyset\). We first compare the weakly
\(r\)-balanced and \(r\)-balanced posets, and then identify the latter with
the proper face poset of a polytope of balancing weights.

\begin{lemma}
\label{weak-balanced-balanced}
Let \(\mathcal B(V,r)=\{S\subsetneq V\mid r\in\conv(S)\}\), ordered by
inclusion, and let
\(\overline{\mathcal B}(V,r)=\{S\subsetneq V\mid r\in\relint\conv(S)\}\).
Then \(\Delta\mathcal B(V,r)\simeq
\Delta\overline{\mathcal B}(V,r)\).
\end{lemma}

\begin{proof}
Let \(\iota:\overline{\mathcal B}(V,r)\hookrightarrow\mathcal B(V,r)\) be the
inclusion. For \(S\in\mathcal B(V,r)\), let \(G_r(S)\) be the minimal face of
\(\conv(S)\) containing \(r\), and put \(\rho(S):=S\cap G_r(S)\).
Since \(G_r(S)=\conv(\rho(S))\) and \(r\in\relint G_r(S)\), we have
\(\rho(S)\in\overline{\mathcal B}(V,r)\). If \(S\subset T\), then
\(G_r(T)\cap\conv(S)\) is a face of \(\conv(S)\) containing \(r\), so
\(G_r(S)\subset G_r(T)\) and \(\rho(S)\subset\rho(T)\). Thus \(\rho\) is
order-preserving. Moreover,
\(\rho\iota=\id_{\overline{\mathcal B}(V,r)}\) and
\(\iota\rho\leq\id_{\mathcal B(V,r)}\).
Pointwise comparable order-preserving maps induce homotopic maps on order
complexes, so \(\rho\) and \(\iota\) induce homotopy inverse maps.
\end{proof}

We now introduce the polytope of weight vectors:
\[
    \Lambda(V,r):=
    \{\lambda\in\mathbb R^m:
    \lambda_i\geq 0,\ 
    \sum_{i=1}^m\lambda_i=1,\ 
    \sum_{i=1}^m\lambda_i v_i=r\}.
\]
Since \(r\in\relint\conv(V)\), this polytope contains a point
\(\mu\) with \(\mu_i>0\) for every \(i\). The affine equations defining
\(\Lambda(V,r)\) have rank \(k+1\), and the strictly positive point shows
that the inequalities \(\lambda_i\geq0\) do not lower its dimension. Hence
\(\dim\Lambda(V,r)=m-k-1\).

For \(\lambda\in \Lambda(V,r)\), put
\(\supp(\lambda)=\{v_i\in V\mid \lambda_i>0\}\).
A subset \(S\subset V\) is weakly \(r\)-balanced if and only if
\(\supp(\lambda)\subset S\) for some \(\lambda\in\Lambda(V,r)\), and it is
\(r\)-balanced if and only if \(\supp(\lambda)=S\) for some such
\(\lambda\). Here we use the standard fact that a point belongs to
\(\relint\conv(S)\) if and only if it has a convex representation with
strictly positive coefficients on all points of \(S\).

For \(S\subset V\), define
\(\Lambda_S:=\{\lambda\in\Lambda(V,r)\mid \supp(\lambda)\subset S\}\).
Equivalently, \(\Lambda_S\) is obtained from \(\Lambda(V,r)\) by imposing the
equations \(\lambda_i=0\) for all \(v_i\notin S\).

\begin{proposition}
\label{balanced-face-polytope}
The map \(S\mapsto\Lambda_S\) is an isomorphism from
\(\overline{\mathcal B}(V,r)\) onto the poset of non-empty proper faces of
\(\Lambda(V,r)\).
\end{proposition}

\begin{proof}
Let \(S\) be \(r\)-balanced, and choose \(\lambda\in\Lambda(V,r)\) with
\(\supp(\lambda)=S\). Then \(\Lambda_S\) is a non-empty face of
\(\Lambda(V,r)\): it is obtained by imposing \(\lambda_i=0\) for
\(v_i\notin S\). Moreover, \(\lambda_i>0\) for every \(v_i\in S\), so
\(\lambda\in\relint\Lambda_S\).

Now let \(G\) be a non-empty face of \(\Lambda(V,r)\), choose
\(\lambda\in\relint G\), and set \(S=\supp(\lambda)\). The equation
\(r=\sum_i\lambda_i v_i\), with \(\lambda_i>0\) precisely for
\(v_i\in S\), shows that \(S\) is \(r\)-balanced. The preceding argument
also gives \(\lambda\in\relint\Lambda_S\). Since the relative interiors of
distinct faces are disjoint, \(G=\Lambda_S\). Thus the inverse map sends a
face to the support of any point in its relative interior. Also,
\(\Lambda_S\subseteq\Lambda_T\) if and only if \(S\subseteq T\).

Finally, \(\Lambda_V=\Lambda(V,r)\). If \(S\subsetneq V\), the point \(\mu\)
with full support does not belong to \(\Lambda_S\), so \(\Lambda_S\) is a
proper face. Hence the restriction to proper balanced subsets gives exactly
the non-empty proper faces.
\end{proof}

\begin{proof}[Proof of Theorem~\ref{ordercomplex}]
Lemma~\ref{weak-balanced-balanced} gives
\(\Delta\mathcal B(V,r)\simeq
\Delta\overline{\mathcal B}(V,r)\).
Since \(\mathcal B(V,r)\neq\emptyset\), the polytope \(\Lambda(V,r)\) has
positive dimension: a proper weakly balanced subset gives a point of
\(\Lambda(V,r)\) with a zero coordinate, distinct from the full-support point
\(\mu\). By Proposition~\ref{balanced-face-polytope},
\(\Delta\overline{\mathcal B}(V,r)\) is the barycentric subdivision of
\(\partial\Lambda(V,r)\). Since \(\dim\Lambda(V,r)=m-k-1\), it follows that
\[
    \Delta\mathcal B(V,r)\simeq
    \Delta\overline{\mathcal B}(V,r)\simeq S^{m-k-2}.
\]
\end{proof}

We finish the section by recording a consequence for Alexander duals of
\(r\)-unbalanced complexes. Recall that the Alexander dual of a simplicial
complex \(\mathcal K\subset 2^V\) is
\(\mathcal K^*=\{\sigma\subset V\mid V\setminus\sigma\notin\mathcal K\}\).

\begin{theorem}
\label{alexander-dual-sphere}
Let \(V=\{v_1,\dots,v_m\}\subset\mathbb R^d\), and let
\(r\in\relint\conv(V)\). Let
\(\mathcal K(V,r)=\{S\subset V\mid r\notin\conv(S)\}\)
be the \(r\)-unbalanced complex. Assume that
\(\mathcal B(V,r)\neq\emptyset\). Put \(k:=\dim\aff(V)\). Then
\[
    \mathcal K(V,r)\simeq S^{k-1}
    \quad\text{and}\quad
    \mathcal K(V,r)^*\simeq S^{m-k-2}.
\]
\end{theorem}

\begin{proof}
The first statement follows from Theorem~\ref{r-unbalanced}, applied inside
the affine hull of \(V\). For the second,
\(\sigma\in\mathcal K(V,r)^*\) if and only if
\(r\in\conv(V\setminus\sigma)\).
Thus \(\sigma\mapsto V\setminus\sigma\) is an order-reversing bijection from
the non-empty face poset of \(\mathcal K(V,r)^*\) to
\(\mathcal B(V,r)\). A poset and its opposite have isomorphic order
complexes, so the barycentric subdivision of \(\mathcal K(V,r)^*\) is
isomorphic to \(\Delta\mathcal B(V,r)\). The result follows from
Theorem~\ref{ordercomplex}.
\end{proof}

\section*{Acknowledgments}

This work was done during the Summer Research Program at MIPT - LIPS-25 (Phystech School of applied mathematics and computer science). The research is supported by the ``Priority 2030'' strategic academic leadership program.

We thank Andrey Ryabichev and Egor Kosolapov for helpful discussions and
comments. We are grateful to Roman Karasev for reading an earlier draft of the
paper and for his remarks.

\section*{Declaration on the Use of AI}
The author used generative AI tools to improve the exposition and clarity of the
manuscript and to help identify gaps in preliminary versions of the proofs. The main
mathematical ideas and the overall structure of the proofs were developed by the
author.

\end{document}